\documentclass{gtart}
\usepackage{amsmath, amssymb, amsthm, graphicx, psfrag} 

\newtheorem{thm}{Theorem}[section]

\newtheorem{lem}[thm]{Lemma}
\newtheorem{prop}[thm]{Proposition}

\theoremstyle{definition}

\numberwithin{equation}{section}

\newcommand{\al}{\alpha}
\newcommand{\be}{\beta}

\newcommand{\x}{\times}

\newcommand{\Z}{\mathbb Z}

\newcommand{\Q}{\mathbb Q}
\newcommand{\R}{\mathbb R}

\newcommand{\del}{\partial}

\DeclareMathOperator{\rot}{rot}

\begin{document}
\mathsurround=1pt 
\title{Classification of tight contact structures on 
small Seifert 3--manifolds with $e_0\geq 0$}

\authors{Paolo Ghiggini, Paolo Lisca, Andr\'{a}s I.~Stipsicz} 

\address{%
\rm Email: \stdspace \tt ghiggini@mail.dm.unipi.it\ \
lisca@dm.unipi.it\\ stipsicz@math-inst.hu}


\begin{abstract}
We classify positive, tight contact structures on closed Seifert
fibered 3--manifolds with base $S^2$, three singular fibers and
$e_0\geq 0$.
\end{abstract}
\maketitle

\section{Introduction}\label{s:intro}
\begin{figure}[ht]
\begin{center}
\psfrag{0}{$0$}
\psfrag{1}{$-\frac{1}{r_1}$}
\psfrag{2}{$-\frac{1}{r_2}$}
\psfrag{3}{$-\frac{1}{r_3}$}
\includegraphics[width=10cm]{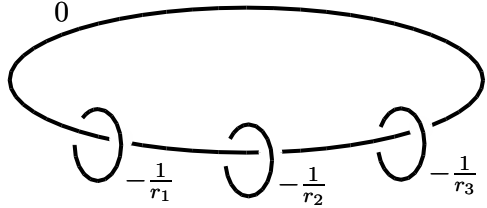}
\end{center}
\caption{Surgery diagram for the Seifert fibered 3--manifold $M(r_1,r_2,r_3)$}
\label{f:seif}
\end{figure}

The classification of positive, tight contact structures on lens
spaces is due to K.~Honda and E.~Giroux~\cite{H1, Gi2}.  Let $M$ be a
closed, small Seifert fibered 3--manifold which is not a lens
space. Then, $M$ has base $S^2$ and exactly three singular
fibers. Equivalently, $M$ is orientation--preserving diffeomorphic
to $M(r_1, r_2, r_3)$ for some $r_1, r_2, r_3\in\Q\setminus\Z$, where
$M(r_1, r_2, r_3)$ denotes the oriented 3--manifold given by the
surgery diagram of Figure~\ref{f:seif}.

Applying Rolfsen twists to the diagram of Figure~\ref{f:seif} it is
easy to show that
\begin{equation}\label{e:twists}
M(r_1,r_2,r_3)=M(r_1+h,r_2+k,r_3-h-k),\quad h,k\in\Z.
\end{equation}
The integer 
\begin{equation*}
e_0(M(r_1,r_2,r_3)):=\sum _{i=1}^3 [r_i]
\end{equation*}
is an invariant of the Seifert fibered 3--manifold $M(r_1, r_2, r_3)$.

 Recently H.~Wu obtained the
classification up to isotopy of positive tight contact structures on
$M(r_1,r_2,r_3)$ (and therefore on small Seifert 3--manifolds)
assuming $e_0\neq -2,-1,0$~\cite{W2}. He used convex surface theory to
derive an upper bound for the number of isotopy classes of tight
contact structures, and produced Legendrian surgery diagrams which
show that the upper bound found in the first step is sharp.

In this note we extend Wu's results to the case $e_0=0$. More
precisely, we classify positive tight contact stuctures on
$M(r_1,r_2,r_3)$ assuming $e_0\geq 0$, by using a set of Legendrian
surgery diagrams which is slightly different from the one used
in~\cite{W2}.

\begin{figure}[ht]
\begin{center}
\psfrag{a}{\small $a^1_0$}
\psfrag{b}{\small $a^2_0$}
\psfrag{c}{\small $a^3_0$}
\psfrag{d}{\small $a^1_1$}
\psfrag{e}{\small $a^2_1$}
\psfrag{f}{\small $a^3_1$}
\psfrag{g}{\small $a^1_2$}
\psfrag{h}{\small $a^2_2$}
\psfrag{i}{\small $a^3_2$}
\psfrag{l}{\small $a^1_{n_1-1}$}
\psfrag{m}{\small $a^2_{n_2-1}$}
\psfrag{n}{\small $a^3_{n_3-1}$}
\psfrag{o}{\small $a^1_{n_1}$}
\psfrag{p}{\small $a^2_{n_2}$}
\psfrag{q}{\small $a^3_{n_3}$}
\psfrag{r}{\small $\cdots$}
\psfrag{1}{\small $K_1$}
\psfrag{2}{\small $K_2$}
\psfrag{3}{\small $K_3$}
\includegraphics[width=12cm]{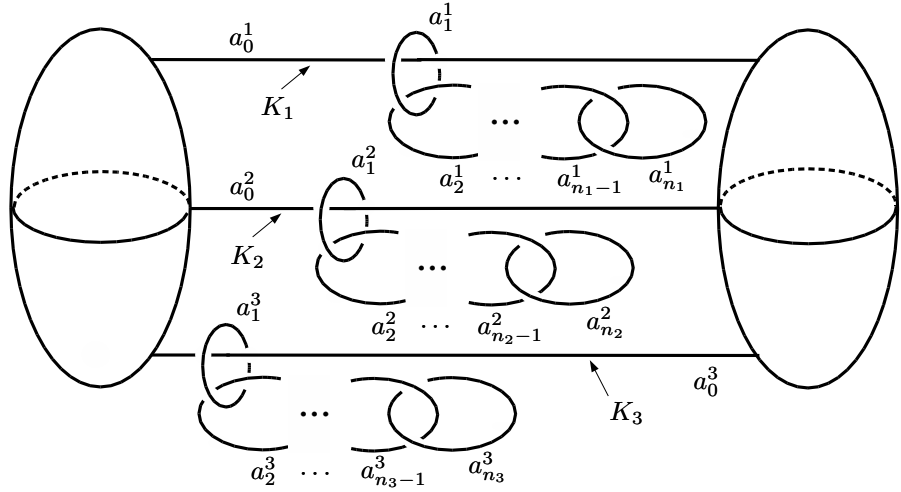}
\end{center}
\caption{Handlebody decomposition of Stein fillings of $M(r_1,r_2,r_3)$}
\label{f:surg}
\end{figure}

Observe that if $e_0(M(s_1,s_2,s_3))\geq 0$, then by~\eqref{e:twists}
we have 
\[
M(s_1,s_2,s_3)=M(r_1,r_2,r_3)\quad\text{for some}\quad r_1 >
0\quad\text{and}\quad 1> r_2, r_3 > 0.
\]
For each of the three rational numbers $r_1$, $r_2$, $r_3$,
we can write
\begin{equation*}
- \frac{1}{r_i}= [a^1_0,a^2_0,\ldots,a^i_{n_i}]:=
a_0^i - \cfrac{1}{a_1^i - \cfrac{1}{\ddots -
  \cfrac{1}{a_{n_i}^i}}},\quad i=1,2,3,
\end{equation*}
for some uniquely determined integer coefficients 
\[
a^1_0\leq -1\quad\text{and}\quad a^2_0, a^3_0,
a^i_1,\cdots,a^i_{n_i}\leq -2.
\]
We define
\[
T(r_1,r_2,r_3):=|(\prod_{i=1}^3 (a_0^i+1)- \prod_{i=1}^3 a_0^i)
\prod_{i=1}^3 \prod_{k=1}^{n_i} (a_k^i+1)|. 
\]
The following is our main result.
\begin{thm}\label{t:main}
Suppose $r_1>0$ and $1>r_2, r_3>0$. Then, $M=M(r_1, r_2, r_3)$ carries
exactly $T(r_1,r_2,r_3)$ positive tight contact structures up to
isotopy. Moreover, each tight contact structure on $M$ has a Stein
filling whose underlying 4--manifold has the handlebody decomposition
given by Figure~\ref{f:surg}.
\end{thm}

\rk{Acknowledgments.} {The first and second authors are members of
  EDGE, Research Training Network HPRN-CT-2000-00101, supported by The
  European Human Potential Programme. The first author was partially
  supported by an EDGE fellowship while he stayed at the Department of
  Mathematics of the Ludwig-Maximilians-Universit{\"a}t Munich. The
  second author was partially supported by MURST. The third author was
  partially supported by OTKA T034885. Part of this work was carried
  out while the third author visited the University of Pisa. He wishes
  to thank the Geometry Group of the Department of Mathematics
  ``L. Tonelli'' for the hospitality.}

\section{Upper bounds}

In this section we establish an upper bound on the number of
isotopy classes of tight contact structures on the Seifert fibered
3--manifold $M=M(r_1,r_2,r_3)$. Let $\xi $ be a tight contact
structure on $M$. Then, a Legendrian knot in $M$ smoothly isotopic to
a regular fiber admits two framings: one coming from the fibration and
the other one coming from the contact structure $\xi$. The difference
between the contact framing and the fibration framing is the
\emph{twisting number} of the Legendrian curve. We say that $\xi$ has
\emph{maximal twisting equal to zero} if there is a Legendrian knot
isotopic to a regular fiber and having twisting number zero.

\begin{prop}[\cite{W1}, Theorem~1.3]\label{p:zerotwisting}
If $r_1, r_2, r_3>0$, then any tight contact structure on
$M(r_1,r_2,r_3)$ has maximal twisting equal to zero.\qed
\end{prop}

We can give an explicit construction of the Seifert manifold
$M(r_1,r_2,r_3)$ as follows. Let $\Sigma$ be an oriented pair of
pants, and identify each connected component of
\[
- \partial (\Sigma \times S^1) = T_1 \cup T_2 \cup T_3 
\]
with $\R^2/\Z^2$, so that $\binom 10$ gives the direction of
$\partial(\Sigma\x\{1\})$ and $\binom 01$ gives the direction of the
$S^1$ factor. Then glue a solid torus $D^2 \times S^1$ to each $T_i$
using the map $\varphi_{A_i}: \partial(D^2 \times S^1) \to T_i$
defined by the matrix
\[
A_i= \left ( \begin{matrix} \alpha_i & \alpha_i' \\ 
- \beta_i & - \beta_i' \end{matrix} \right),
\]
where 
\[
\frac{\beta_i}{\alpha_i}=r_i,\quad 
\alpha_i' \beta_i - \alpha_i \beta'_i=1,\quad\text{and}\quad 
0 < \alpha_i' < \alpha_i. 
\]
Since $r_i > 0$, it follows that $\beta_i>0$. The singular fibers of
the Seifert fibration will be denoted by $F_i$, $i=1,2,3$.

\begin{lem}\label{l:background}
Let $r_1, r_2, r_3>0$, and let $\xi$ be a tight contact structure on
$M=M(r_1,r_2,r_3)$. Then, the singular fibers $F_i$ can be made
Legendrian with twisting number $-1$. Moreover, there exist convex
neighbourhoods $U_i$ of $F_i$ such that $-\del(M\setminus U_i)$ has
infinite slope.
\end{lem}

\begin{proof}
Using Proposition~\ref{p:zerotwisting}, one can isotope $\xi$ until
there is a Legendrian regular fiber $L$ with twisting number
$0$. Then, one can make the singular fibers $F_i$ Legendrian with very
low twisting numbers $n_i<0$. Let $V_i$ be a standard convex
neighbourhood of $F_i$, and let $A_i$ be a convex vertical annulus
between $L$ and a ruling of $\partial (M \setminus V_i)$. By the
Imbalance Principle~\cite[Proposition~3.17]{H1}, $A_i$ produces a
bypass attached to $\partial V_i$ along a Legendrian curve with slope
\[
-\frac{\alpha_i}{\alpha'_i}<-1. 
\]
Using the Twisting Number Lemma~\cite[Lemma~4.4]{H1} we can increase
the twisting number of $F_i$ up to $-1$. Then, we can thicken $V_i$
further in order to obtain a convex solid torus $U_i$ such that
$-\del(M\setminus U_i)$ has infinite slope.
\end{proof}

Let $\Sigma$ be a pair of pants. We say that a tight contact structure
$\xi$ on $\Sigma \times S^1$ is {\em admissible} if there is no
contact embedding $(T^2 \times I, \xi_{\pi}) \hookrightarrow (\Sigma
\times S^1, \xi)$, where $\xi_{\pi}$ is a tight contact structure with
convex boundary and twisting $\pi$ (see~\cite[{\S}~2.2.1]{H1} for the
definition of twisting).

\begin{lem}
Let $\xi$ be a tight contact structure on $M=M(r_1,r_2,r_3)$, and
suppose that the singular fibers $F_i$ are all Legendrian with
twisting number $-1$.  Let $V_i$ be a standard neighbourhood of
$F_i$, for $i=1,2,3$. Then, the restriction of $\xi$ to
$M\setminus(V_1 \cup V_2 \cup V_3)$ is admissible.
\end{lem}

\begin{proof}
Arguing by contradiction, suppose that 
\[
(T^2 \times I, \xi_{\pi})\subset (M\setminus(V_1 \cup V_2 \cup V_3),
\xi).
\]
Any embedded torus in $M \setminus (V_1 \cup V_2 \cup V_3)$ contains a
homotopically nontrivial embedded circle $C$ which bounds a disc $D$
in $M$. Since $(T^2 \times I, \xi_{\pi})$ has twisting $\pi$,
by~\cite[Proposition~4.16]{H1} there is a standard torus $T\subset T^2
\times [0,1]$ with Legendrian divides isotopic to $C$. Then, $D$ is
isotopic to an overtwisted disc in $(M, \xi)$.
\end{proof}

\begin{lem}\label{l:scambio}
Let $\Sigma$ be a pair of pants and $\xi_+$, $\xi_-$ admissible tight
contact structures on $\Sigma \times S^1$. Suppose that the boundary
\[
- \partial (\Sigma \times S^1)=T_1 \cup T_2 \cup T_3 
\]
consists of tori in standard form, with $\#\Gamma_{T_i}=2$, $i=1,2,3$,
and slopes $s_1=s_2=s_3=-1$.  Let $\Sigma'\subset \Sigma$ be another
pair of pants such that
\[
\Sigma \times S^1 = \Sigma' \times S^1 \cup (T_1 \times
I) \cup (T_2 \times I) \cup (T_3 \times I), 
\]
where $\xi_+|_{T_i \times I}$ and $\xi_-|_{T_i \times I}$ are,
respectively, a positive and a negative basic slice with boundary
slopes $-1$ and $\infty$. Then, $(\Sigma \times S^1, \xi_+)$ is
isotopic to $(\Sigma \times S^1, \xi_-)$.
\end{lem}

\begin{proof} 

By~\cite[Lemma~4.13]{GhS} we can find  vertical
annuli $A_{\pm}\subset \Sigma \times S^1$ between $T_1$ and $T_2$ such that
\begin{enumerate}
\item $A_{\pm}$ is convex and has Legendrian boundary with respect to
$\xi_{\pm}$
\item the dividing set of $A_{\pm}$ has no boundary parallel curves.
\end{enumerate} 
In fact, it is easy to check that $A_-$ and $A_+$ are isotopic. Let
$\phi_t$ be an isotopy of $\Sigma \times S^1$ which is the identity on
the boundary and such that $\phi_1(A_-)=A_+$. To prove the lemma it
suffices to show that $\xi_+$ is isotopic to
$(\phi_1)_*(\xi_-)$. Therefore, without loss of generality we may
assume $A_+=A_-=:A$.

After rounding the edges,
\[
(\Sigma \times S^1 \setminus A, 
\xi_{\pm})
\]
is isomorphic to a tight solid torus $(T^2 \times [0,2], \eta_{\pm})$,
in such a way that $T^2 \times [0,1]$ corresponds to 
\[
\Sigma\x S^1 \setminus (T_3\x I\cup A) 
\]
and $T^2 \times [1,2]$ corresponds to $T_3\x I$.  With this
identification, the slopes of $T^2\x\{0\}$ and $T^2\x\{1\}$ are,
respectively, $-1$ and $+1$. Also, notice that $(T^2 \times [0,2],
\eta_{\pm})$ is minimally twisting because $(\Sigma \times S^1,
\xi_{\pm})$ is admissible.

A relative Euler class computation as in the proof
of~\cite[Lemma~4.13]{GhS} shows that 
\[
(T^2 \times [0,1], \eta_{\pm}|_{T^2 \times [0,1]})\quad\text{and}\quad
(T_1 \times [0,1], \xi_{\pm}|_{T_1 \times [0,1]}) 
\]
are basic slices with the same sign. Similarly, 
\[
(T^2 \times [1,2], \eta_{\pm}|_{T^2 \times [1,2]})\quad\text{and}\quad
(T_3 \times [0,1], \xi_{\pm}|_{T_3 \times [0,1]}) 
\]
have opposite signs, because the diffeomorphism between $T_3\x\{0\}$
and $T^2 \times \{ 2 \}$ reverses orientations. Thus, $(T^2 \times
[0,2], \eta_+)$ and $(T^2 \times [0,2], \eta_-)$ both decompose into a
positive basic slice and a negative basic slice belonging to the same
continued fraction block, and therefore by~\cite[Lemma~4.14]{H1} they
are isotopic. We conclude that $(\Sigma \times S^1, \xi_+)$ and
$(\Sigma \times S^1, \xi_-)$ are isotopic because they decompose into
pairwise isotopic pieces.
\end{proof}

\begin{lem}\label{l:numeri}
The following equalities hold for every $i=1,2,3$:
\begin{equation}\label{e:firsteq}
-\frac{\alpha_i}{\alpha'_i}=[a^i_{n_i},\ldots, a^i_0]
\end{equation}
\begin{equation}\label{e:secondeq}
-\frac{\alpha_i +(a^i_0+1)\beta_i}{\alpha'_i +(a^i_0+1)\beta'_i}=
[a^i_{n_i},\ldots, a^i_1+1]
\end{equation}
\end{lem}

\begin{proof}
Equality~\eqref{e:firsteq} follows from~\cite[Lemma~A4]{OW}. A
straightforward computation gives
\[
-\frac{\al_i+(a^i_0+1)\be_i}{-\al_i-a^i_0\be_i} = 
[a^i_1+1,a^i_2,\ldots,a^i_{n_i}].
\]
Thus,~\cite[Lemma~A4]{OW} together with the fact that
\[
(\al'_i+(a^i_0+1)\be'_i)(-\al_i-a^i_0\be_i) - 
(\al_i+(a^i_0+1)\be_i)(-\al'_i - a^i_0\be'_i) = 1 
\]
implies Equality~\eqref{e:secondeq}.
\end{proof}

\begin{lem}\label{l:firstborder}
Suppose $a^1_0<-1$. Then, for every $i=1,2,3$, the slope of the border
between the first (i.e.~the outermost) and the second continued
fraction block of $\partial U_i$, computed in the basis of $-\partial
(M \setminus U_i)$, is
\[
\frac{1}{a_0^i+1}.
\]
\end{lem}

\begin{proof}
According to~\cite[{\S}~4.4.4]{H1} and in view of
Equation~\eqref{e:secondeq}, the slope of the border between the first
and the second continued fraction block of $\partial U_i$ is
\[
[a_{n_i}^i, \ldots ,a_1^i+1] = 
- \frac{\alpha_i +(a_0^i+1) \beta_i}{\alpha_i' +(a_0^i+1) \beta_i}. 
\]
A direct computation via the matrix $A_i$ gives the slope in the
basis of $-\partial (M \setminus U_i)$.
\end{proof}

\begin{thm}\label{t:upperbound}
Suppose $r_1>0$ and $1>r_2, r_3>0$. Then, $M(r_1, r_2, r_3)$ 
carries at most
\[
T(r_1,r_2,r_3):=|(\prod_{i=1}^3 (a_0^i+1)- \prod_{i=1}^3 a_0^i)
\prod_{i=1}^3 \prod_{k=1}^{n_i} (a_k^i+1)|
\]
distinct tight contact structures up to isotopy.
\end{thm}

\begin{proof}
Fix a decomposition of $M= M(r_1, r_2, r_3)$ as in
Lemma~\ref{l:background}:
\[
M=M \setminus (U_1 \cup U_2 \cup U_3)\cup_{i=1}^3 U_i.
\]
Let $V_i\subset U_i$ be a standard neighborhood of the singular fiber
$F_i$, for $i=1,2,3$. Then, up to an isotopy which is fixed on the
boundary there is at most one tight contact structure on $\cup_{i=1}^3
V_i$. Moreover, by~\cite[Lemma~11]{EH}, there is at most one
admissible tight contact structure on $M \setminus (U_1 \cup U_2 \cup
U_3)$ up to an isotopy not necessarily fixed on the boundary. Notice
that, in general, one should allow only isotopies which are fixed on
the boundary. But in our situation there is no loss in allowing more
general isotopies on $M \setminus (U_1 \cup U_2 \cup U_3)$ because,
as~\cite[Lemma~4.4.]{Gh1} shows, any isotopy on $\del U_i$ extends to
$U_i$.

Let $N_k^i$ be the $(k+1)$-th continued fraction block of $(U_i,
\xi|_{U_i})$. The proof of Lemma~\ref{l:background} shows that the
solid tori $U_i$ have boundary slope
\[
-\frac{\alpha_i}{\alpha_i'}< -1. 
\]
Therefore, by Equation~\eqref{e:firsteq} and~\cite[{\S}~4.4.4]{H1} the slope
of the border between $N_{k-1}^i$ and $N_k^i$ is
$s_k^i=[a^i_{n_i},\ldots, a^i_k+1]$. Thus, if we define
\[
p_k^i := \# \{ \hbox{positive basic slices in } N_k^i \},
\]
we have
\begin{equation}\label{e:pineq}
0 \leq p_0^i \leq |a_0^i+1|\quad\text{and}\quad 0\leq p_k^i \leq
|a_k^i+2|.
\end{equation}

Let $V'_i:=U_i\setminus N_0^i$. By Inequalities~\eqref{e:pineq}
and~\cite[Theorem~2.2]{H1}, there are exactly 
\[
|(a_1^i+1) \cdots (a_n^i+1)|
\]
distinct isotopy classes of tight contact structures on $V_i'
\setminus V_i$, and $|a_0^1 a_0^2 a_0^3|$ possible configurations of
signs $(p_0^1, p_0^2, p_0^3)$ in $M \setminus (V_1' \cup V_2' \cup
V_3')$. This immediately gives the number
\begin{equation}\label{e:upperbound}
|a_0^1 a^2_0 a^3_0\prod_{i=1}^3\prod_{k=1}^{n_1} (a^i_k+1)|
\end{equation}
as an upper bound for the number of isotopy classes of tight contact
structures on $M$. If $a_0^1=-1$ (which is equivalent to $e_0(M)>0$),
clearly the quantity in~\eqref{e:upperbound} coincides with
$T(r_1,r_2,r_3)$, and the statement follows.

Thus, we are left to prove the statement when $a_0^1<-1$ (which is
equivalent to $e_0(M)=0$). In this case, the upper bound given
by~\eqref{e:upperbound} is not optimal, because Lemma~\ref{l:scambio}
shows that different sign configurations do not necessarily yield
distinct contact structures on $M \setminus (V_1' \cup V_2' \cup
V_3')$. In fact, if $N_0^1$, $N_0^2$ and $N_0^3$ contain basic slices
$B_1$, $B_2$ and $B_3$ with the same sign, by \cite[Lemma~4.14]{H1} we
can arrange the basic slice decomposition of each $N_0^i$ so that
$B_i$ is the first basic slice. Thus, it is easy to check using
Lemma~\ref{l:firstborder} that each $B_i$ has boundary slopes $-1$ and
$\infty$ when computed in the basis of $-\partial (M \setminus
U_i)$. Applying Lemma~\ref{l:scambio}, we are allowed to change the
sign of all three basic slices simultaneously without changing the
isotopy type of the contact structure. This shows that the
configuration $(p_0^1, p_0^2, p_0^3)$ is equivalent to $(p_0^1 \pm 1,
p_0^2 \pm 1, p_0^3 \pm1)$ (with the same signs chosen in each slot),
whenever the sums are defined.

We can easily count the different possibilities for $p_0^i$: by the
above argument we can always arrange that one of the $p_0^i$'s is
maximal, i.e.~equal to $\vert a_0^i+1\vert$. For the other two we have
$\vert a_0^j\vert \cdot \vert a_0^k\vert $ many choices (where $\{
i,j,k\} =\{ 1,2,3\}$). A simple computation shows that the total
number of possibilities is equal to
\[
\vert a_0^1\vert \cdot \vert a_0^2\vert  +
\vert a_0^2\vert \cdot \vert a_0^3\vert +
\vert a_0^3\vert \cdot \vert a_0^1\vert -
\vert a_0^1\vert -  \vert a_0^2\vert - \vert a_0^3\vert +1,
\]
and this expression is equal to
\[
\vert (\prod_{i=1}^3 (a_0^i+1)- \prod_{i=1}^3 a_0^i)\vert.
\]
This proves the statement when $a_0^1<-1$, and concludes the proof.
\end{proof}   

\section{Lower bounds}

In this section we construct $T(r_1,r_2,r_3)$ distinct isotopy classes
of Stein fillable, hence tight, contact structures on
$M=M(r_1,r_2,r_3)$ assuming $r_1, r_2, r_3>0$.

Notice that the diagram of Figure~\ref{f:surg} gives the
handlebody decomposition of a 4-manifold $X$ with boundary
diffeomorphic to $M$. The decomposition involves a single 1--handle
and some 2--handles. 

Since $a_0^i\leq -1$ and $a_k^i\leq -2$ for $k>0$, following~\cite{G}
it is easy to describe a Stein structure on $X$ by putting the knots
into Legendrian position and stabilizing them until the prescribed
framing coefficient becomes $-1$ with respect to their contact
framing. This way we get
\[
| \prod_{i=1}^3 \prod_{k=0}^{n_i} (a_k^i+1)| 
\]
different Legendrian diagrams giving Stein structures on $X$ and
therefore tight contact structures on $\partial X=M$. Let us denote by
$\xi_J$ the contact structure corresponding to a Stein structure $J$
on $X$. According to~\cite{LM}, if $c_1(J_1)\neq c_1(J_2)$ then the
induced contact structures $\xi_{J_1}$ and $\xi_{J_2}$ are
nonisotopic. Our aim is to count the number of distinct first Chern
classes obtainable in this way.

In order to do this, we start by fixing a basis of $H_2(X;\Z)$.  We
will present the second homology group using cellular homology. It is
well--known~\cite{GS} that the framed knots in the diagram correspond
to 2--cells. Hence, a choice of orientation for the knots gives rise
to a basis for the group $C_2(X)$ of 2--chains for $X$.

Let $C_1(X)$ denote the group generated by the 1--cells, i.e.~by the
1--handles in the handle decomposition.  Since there are no 3--handles
present in the handle decomposition, $H_2(X;\Z)$ can be computed as
the kernel of the map $\varphi\colon C_2(X)\to C_1(X)$, given on a basis
element $K\in C_2(X)$ corresponding to the knot $K$ as
\[
\varphi (K)=\sum a_i L_i, 
\]
where $L_i$ runs through all 1--handles and $a_i\in \Z$ is the
algebraic number of times $K$ passes through the 1--handle $L_i$. In
our case we have $C_1(X)\cong\Z$ and, for a suitable choice of
orientations,
\[
\varphi(K_i)=1\quad\text{for}\quad i=1,2,3, 
\]
where the knots $K_i$ are indicated in Figure~\ref{f:surg}, and
$\varphi$ is zero on each basis element $K\neq K_1, K_2,
K_3$. Consequently, a basis of $H_2(X; \Z)$ can be given by the
homology classes corresponding to the unknots of Figure~\ref{f:surg}
together with the classes of
\[
K_1-K_2, K_1-K_3\in C_2(X).
\]
It follows from the results of~\cite{G} that if $K\neq K_1, K_2,
K_3$, then
\[
\langle c_1(J), [K]\rangle = \rot(K), 
\]
while if $\{i,j\}=\{1,2\}$ or $\{i,j\}=\{1,3\}$,
\[
\langle c_1(J), [K_i-K_j]\rangle = 
\rot(K_i)-\rot(K_j).
\]
Theorem~\ref{t:main} follows immediately from
Theorem~\ref{t:upperbound} together with the following

\begin{prop}\label{p:lowerbound}
Suppose $r_1, r_2, r_3>0$. Then, $M(r_1, r_2, r_3)$ carries at least 
\[
T(r_1,r_2,r_3):=
|(\prod_{i=1}^3 (a_0^i+1)- \prod_{i=1}^3 a_0^i)
\prod_{i=1}^3 \prod_{k=1}^{n_i} (a_k^i+1)| 
\]
distinct Stein fillable contact structures up to isotopy. 
\end{prop}

\begin{proof}
Let $J_1$ and $J_2$ be Stein structures on $X$ resulting from oriented
Legendrian surgery diagrams as above. Denote by $r_i^k(1)$
and $r_i^k(2)$ ($i=1,2,3$, $k=0,\ldots, n_i$) the rotation numbers 
of the Legendrian knots appearing in the two diagrams. It follows from
the above discussion that if either 
\[
r^0_1(1)-r^0_2(1)\neq r^0_1(2)-r^0_2(2),\quad
r^0_1(1)-r^0_3(1)\neq r^0_1(2)-r^0_3(2), 
\]
or
\[
r_i^k(1)\neq r_i ^k(2)\quad\text{for some}\quad k>0, 
\] 
then $c_1(J_1)\neq c_1(J_2)$, and therefore by~\cite{LM} the induced
contact structures $\xi_1$ and $\xi_2$ on $\partial X=M$ are not
isotopic. The conclusion follows from a computation similar to the one
given in the proof of Theorem~\ref{t:upperbound}.
\end{proof}

\end{document}